\documentclass[12pt,a4paper]{amsart}
\usepackage{amssymb,mathrsfs,amscd}

\let\g=\mathfrak

\newcommand{\0}{{\boldsymbol{0}}}
\newcommand{\1}{{\boldsymbol{1}}}

\newcommand{\Ad}{\operatorname{Ad}}

\newcommand{\Spf}{\operatorname{Spf}}

\newcommand{\tr}{\operatorname{tr}}

\newcommand{\Vect}{\operatorname{Vect}}
\newcommand{\End}{\operatorname{End}}
\newcommand{\defn}{\mathop{=}\limits^{\text{def}}}

\def\BC{{\mathcal B}}
\def\CC{{\mathcal C}}
\def\DC{{\mathcal D}}

\def\LC{{\mathcal L}}

\def\NC{{\mathcal N}}

\def\SC{{\mathcal S}}

\def\UC{{\mathcal U}}
\def\VC{{\mathcal V}}

\let\a=\alpha

\let\l=\lambda
\let\L=\Lambda

\let\W=\Omega
\let\r=\rho

\let\d=\delta
\let\D=\Delta
\let\de=
\partial

\newtheorem{theo}{Theorem}[section]
\newtheorem{defi}{Definition}[section]

\newtheorem{lemme}{Lemma}[section]

\let\what=\widehat
\def\som{\mathop{\sum}\limits}

\def\tens{\mathop{\otimes}\limits}
\def\sdir{\mathop{\oplus}\limits}
\def\mult{\mathop{\times}\limits}

\let\oo=\infty

\def\hfl#1#2{\smash{\mathop{\hbox to
12mm{\rightarrowfill}}\limits^{\scriptstyle#1}_{\scriptstyle#2}}}

\def\ihfl#1#2{\smash{\mathop{\hbox to
12mm{\hookrightarrowfill}}\limits^{\scriptstyle#1}_{\scriptstyle#2}}}

\def\lihfl#1#2{\smash{\mathop{\hbox to
12mm{\hookleftarrowfill}}\limits^{\scriptstyle#1}_{\scriptstyle#2}}}

\def\lhfl#1#2{\smash{\mathop{\hbox to
12mm{\leqslantftarrowfill}}\limits^{\scriptstyle#1}_{\scriptstyle#2}}}

\title[Invariant generalized functions]{Invariant generalized functions on $\g{sl}(2,\mathbb R)$ with values in a $\g{sl}(2,\mathbb R)$-module}
\author{P. Lavaud}
\date{\today}
\begin{document}

\maketitle

\begin{abstract}
Let $\g g$ be a finite dimensional real Lie algebra. Let $\r:\g g\to\End(V)$ be a representation of $\g g$
in a finite dimensional real vector space. 
Let $\CC_{V}=\big(\End(V)\tens S(\g g)\big)^{\g g}$ be the algebra of $\End(V)$-valued invariant differential operators with constant coefficients on $\g g$. Let $\UC$ be an  open subset of $\g g$. 
  We consider the problem of determining the space   of generalized functions  $\phi$ on $\UC$ with values in $V$ which are locally  invariant and such that $\CC_{V}\phi$ is finite dimensional.
 
In this article we consider the case $\g g=\g{sl}(2,\mathbb{R})$. Let $\NC$ be the nilpotent cone of $\g{sl}(2,\mathbb{R})$. We prove that when $\UC$ is $SL(2,\mathbb {R})$-invariant, then $\phi$ is determined by its restriction to $\UC\setminus\NC$  where $\phi$ is analytic (cf. Theorem \ref{th:Principal}). In general this is false when $\UC$ is not $SL(2,\mathbb {R})$-invariant and $V$ is not trivial. Moreover, when $V$ is not trivial, $\phi$ is not always locally $L^{1}$. Thus, this case is different and more complicated than the situation considered by Harish-Chandra (cf. \cite{HC64b,HC65a}) where $\g g$ is reductive and $V$ is trivial.

To solve this problem we find all the locally invariant generalized functions supported in the nilpotent cone $\NC$. We do this locally in a neighborhood   of  a nilpotent element $Z$ of $\g g$ (cf. Theorem \ref{th:local}) and
on an $SL(2,\mathbb {R})$-invariant open subset $\UC\subset\g{sl}(2,\mathbb{R})$  (cf. Theorem \ref{th:global}). Finally, we also give an application of our main theorem to the Superpfaffian (cf. \cite{Lav04}).

\end{abstract}

\section{Introduction}

Let $\g g$ be a finite dimensional real Lie algebra. Let $\r:\g g\to\End(V)$ be a representation of $\g g$
in a finite dimensional real vector space. 
Let $\CC_{V}=\big(\End(V)\tens S(\g g)\big)^{\g g}$ be the algebra of $\End(V)$-valued invariant differential operators  with constant coefficients on $\g g$. It is the \emph{classical family algebra} in the terminology of Kirillov (cf. \cite{Kir00}). Let $\UC$ be an  open subset of $\g g$. 
  We consider the problem of determining  the space   of generalized functions $\phi$ on $\UC$ with values in $V$ which are locally invariant and such that $\CC_{V}\phi$ is finite dimensional.
  
When $V=\mathbb {R}$ is the trivial module and $\g g$ is reductive, the problem was solved by Harish-Chandra (cf. in  particular \cite{HC64b,HC65a}). Let $\phi$ be a locally invariant generalized function  such that $S(\g g)^{\g g}\phi$ is finite dimensional. He proved that $\phi$ is locally $L^{1}$, $\phi$ is determined by its restriction $\phi|_{\g g'}$ to the open subset $\g g'$ of semi-simple regular elements of $\g g$ and $\phi|_{\g g'}$ is analytic.

\medskip

In this article we consider the case $\g g=\g{sl}(2,\mathbb{R})$. Let $\NC$ be the nilpotent cone of $\g{sl}(2,\mathbb{R})$. In this case $\g g'=\g{sl}(2,\mathbb{R})\setminus \NC$. Let $\phi$ be a locally invariant generalized function  on $\UC$ with values in $V$ such that $\CC_{V}\phi$ is finite dimensional. We prove that when $\UC$ is $SL(2,\mathbb {R})$-invariant, then $\phi$ is determined by its restriction to $\UC\setminus\NC$  where $\phi$ is analytic (cf. Theorem \ref{th:Principal}). In general this is false when $\UC$ is not $SL(2,\mathbb {R})$-invariant and $V$ is not trivial. Moreover, when $V$ is not trivial, $\phi$ is not always locally $L^{1}$. Finally, we also give an application of our main theorem to the Superpfaffian (cf. \cite{Lav04}).

\bigskip

To solve the problem we find all the locally invariant generalized functions supported in the nilpotent cone $\NC$. Let $V_{n}$ be the $n+1$-dimensional irreducible representation of $\g{sl}(2,\mathbb{R})$.  Let $\UC$ be an open subset of $\g{sl}(2,\mathbb{R})$.  We denote by $\CC^{-\oo}(\UC,V_{n})^{\g{sl}(2,\mathbb{R})}$ the set of locally invariant generalized functions on $\UC$ with values in $V_{n}$.
Let $\square$ be the Casimir operator on $\g g$.
 
 \bigskip
 
We denote by $\NC^+$ (resp. $\NC^{-}$) the ``upper'' (resp. ``lower'') half nilpotent cone (cf. \ref{sec:Notations}). We put:
\begin{eqnarray}
\SC_{n}^{0}(\UC)&=&\{\phi\in\CC^{-\oo}(\UC,V_{n})^{\g{sl}(2,\mathbb{R})}\,/\, \phi|_{\UC\setminus\{0\}}=0\};\\
\SC_{n}^{\pm}(\UC)&=&\{\phi\in\CC^{-\oo}(\UC,V_{n})^{\g{sl}(2,\mathbb{R})}\,/\, \phi|_{\UC\setminus(\NC^{\pm}\cup\{0\})}=0\};\\
\SC_{n}(\UC)&=&\{\phi\in\CC^{-\oo}(\UC,V_{n})^{\g{sl}(2,\mathbb{R})}\,/\, \phi|_{\UC\setminus\NC}=0\}.
\end{eqnarray}

\bigskip

 Let $Z\in\NC^+$. We assume that  $\UC$ is a suitable open neighborhood  of $Z$ (cf. section \ref{sec:s_n}).
  Let $\d_{\NC^\pm}$ be an invariant generalized function  with support $\NC^\pm\cup\{0\}$ (cf. section \ref{sec:Dirac}).
We construct an invariant function $s_{n}$ on $\NC\cap\UC$ with values in $V_{n}$.
 We prove (cf. Theorem \ref{th:local}):
\begin{enumerate}
\item When $n$ is even, $\SC_{n}(\UC)$ is an infinite dimensional vector space  with basis:
\begin{equation}
\big(\square^k(s_{n}\d_{\NC^+})\big)_{k\in\mathbb{N}}.
\end{equation}

\item When $n$ is odd, $\dim(\SC_{n}(\UC))=\frac{n+1}{2}$ and a basis is given by:
\begin{equation}
\big(\square^k(s_{n}\d_{\NC^+})\big)_{0\le k\le \frac{n-1}{2}}.
\end{equation}

\end{enumerate}

\bigskip

We assume that $\UC$ is an $SL(2,\mathbb {R})$-invariant open subset of $\g{sl}(2,\mathbb{R})$. If  $\UC\cap \NC\not=\emptyset$, we have $\NC^+\subset\UC$ or $\NC^{-}\subset \UC$. 
We prove (cf. Theorem \ref{th:global}):

\begin{enumerate}
\item 
\begin{equation}
\begin{cases}
\SC_{n}^{0}(\UC)=\{0\}&\text {if }0\not\in\UC;\\
\SC_{n}^{0}(\UC)\simeq (V_{n}\tens S(\g{sl}(2,\mathbb{R})))^{\g{sl}(2,\mathbb{R})}
&\text{if }0\in\UC.
\end{cases}
\end{equation}

\item When $n$ is even, we have:
\begin{eqnarray}
\SC_{n}(\UC)&=&\SC_{n}^{0}(\UC)\oplus\Vect\{\square^k(s_{n}\d_{\NC^{+}})|_{\UC}/k\in\mathbb{N}\}
\oplus\Vect\{\square^k( s_{n} \d_{\NC^{-}})|_{\UC}/k\in\mathbb{N}\}
\\
\SC_{n}^{\pm}(\UC)&=&\SC_{n}^{0}(\UC)\oplus\Vect\{\square^k(s_{n}\d_{\NC^{\pm}})|_{\UC}/k\in\mathbb{N}\}
\end{eqnarray}

\item When $n$ is odd:
\begin{equation}
\SC_{n}(\UC)=\SC_{n}^{\pm}(\UC)=\SC_{n}^{0}(\UC)
\end{equation}

\end{enumerate}

\bigskip

Finally, let $\UC$ be an  open subset of $\g{sl}(2,\mathbb{R})$. Let $V$ be the space of a real finite dimensional representation of $\g g$.
Let $\phi$ be an invariant function defined on $\UC$  such that 
$\CC_{V}\phi$ is finite dimensional. This last condition is equivalent to the existence of $r\in\mathbb {N}$ and $(a_{0},\dots,a_{r-1})\in\mathbb {R}^r$ such that:
\begin{equation*}
\Big(\square^r+\som_{k=0}^{r-1}a_{k}\square^k\Big)\phi=0.
\end{equation*}
Moreover, we assume  that
$\phi|_{\UC\setminus\NC}=0$.
 We prove (cf. Theorem \ref{th:InvGlobal}) that
if $\UC$ is $SL(2,\mathbb {R})$-invariant, then we have $\phi=0$.

\smallskip

In general, when $\UC$ is not $SL(2,\mathbb {R})$-invariant, there exist non trivial solutions of the equation  $\square^k\phi=0$ which are supported in the nilpotent cone (cf. Theorem \ref{th:InvLocal}).

\bigskip

{\bf Acknowledgment:}   I wish to thank Michel Duflo for many fruitful discussions on the subject  and many useful comments  on preliminary versions of this article.

\tableofcontents

\section{Notations}

Let $\g g$ be a finite dimensional real Lie algebra. Let $\r:\g g\to\End(V)$ be a representation of $\g g$
in a finite dimensional real vector space $V$.
Let $\UC$ be an open subset of $\g g$.
We denote by $\DC_{c}^{\oo}(\UC)$ the space of compactly supported  smooth densities on $\UC$. We put:
\begin{equation}
\CC^{-\oo}(\UC,V)=\LC(\DC_{c}^{\oo}(\UC),V),
\end{equation}
where $\LC$ stands for continuous homomorphisms. It is the space of generalized functions on $\UC$ with values in $V$. We put $\CC^{- \oo}(\UC)=\CC^{- \oo}(\UC,\mathbb {R})$. For $\phi\in\CC^{-\oo}(\UC,V)$ and $\mu\in\DC_{c}^{\oo}(\UC)$, we denote by:
\begin{equation}
\int_{\UC}\phi(Z)d\mu(Z)
\end{equation}
the image of $\mu$ by $\phi$. We have:
\begin{equation}
\CC^{-\oo}(\UC,V)= \CC^{-\oo}(\UC)\tens V
\end{equation}
(we will also write $\phi v$ for $\phi\tens v$).
\bigskip

Let $Z\in\g g$. We denote by $\de_{Z}$ the derivative in the direction $Z$. It acts on $\CC^{-\oo}(\UC)$ and on  $\CC^{-\oo}(\UC,V)$. We extend $\de$ to a morphism  of algebras from $S(\g g)$ to the algebra of differential operators with constant coefficients on $\g g$.
We denote by $\LC_{Z}$ the differential operator defined by:  
\begin{equation}
(\LC_{Z}\phi)(X)=\frac{d}{dt}\phi(X-t[Z,X])\big|_{t=0}.
\end{equation}
The map $Z\mapsto \LC_{Z}$ is a Lie algebra homomorphism from $\g g$ into  the algebra of differential operators on $\g g$. 
Let $Z\in\g g$ and $\phi\tens v\in\CC^{-\oo}(\UC)\tens V$ , we put:
\begin{equation}
Z.( \phi\tens v)= \phi\tens \r(Z)v+ (\LC_{Z}\phi)\tens v.
\end{equation}
 In other words, if we extend $\LC_{Z}$ (resp. $\r(Z)$) linearly to a representation  of $\g g$ in $\CC^{-\oo}(\UC,V)$, we have for $\phi\in\CC^{-\oo}(\UC,V)$:
\begin{equation}
Z.\phi=(\r(Z)+\LC_{Z})\phi.
\end{equation}

\bigskip

We say that $\phi\in\CC^{-\oo}(\UC,V)$ is locally invariant if 
 for any $ Z\in\g g$ we have $ Z.\phi=0.$
We put:
\begin{equation}
\CC^{-\oo}(\UC,V)^{\g g}=\{\phi\in\CC^{-\oo}(\UC,V)\,/\,\forall Z\in\g g,\,Z.\phi=0\}.
\end{equation}

\section{Support $\{0\}$ distributions}

\medskip

In this section  we assume that $\g g$ is unimodular.
We choose an invariant measure $dZ$ on $\g g$. We define the Dirac function $\d_{0}$  on $\g g$ with support $\{0\}$ (which depends on the choice of $dZ$) by the following.
Let $\CC_{c}^{\oo}(\g g)$ be the set of smooth compactly supported functions on $\g g$. Then:
\begin{equation}
\forall f\in\CC^{\oo}_{c}(\g g),\int_{\g g}\d_{0}(Z)f(Z)dZ=f(0).
\end{equation}

\medskip

We have the following well known theorem:

\begin{theo}\label{th:Supp0}
Let $\g g$ be a finite dimensional unimodular real Lie algebra and $V$ be a finite dimensional $\g g$-module. Then:

\begin{equation}
\big\{\phi\in\CC^{-\oo}(\g g,V)^{\g g}/\phi|_{\g g\setminus \{0\}}=0\}\simeq (V\tens S(\g g))^{\g g}.
\end{equation}
The isomorphism (which depends on the choice of $dZ$) sends $\som_{i}v_{i}\tens D_{i}\in  (V\tens S(\g g))^{\g g}$ to
$\som_{i}\big(\de_{D_{i}}\d_{0}\big)v_{i}.$

\end{theo}

\bigskip

\section{Support in the nilpotent cone}

From now on,  we assume that $\g g=\g{sl}(2,\mathbb{R})$. 

\subsection{Notations}\label{sec:Notations}

We put:
\begin{equation}
H=\left(\begin{array}{cc}1 & 0 \\0 & -1\end{array}\right);\qquad
X=\left(\begin{array}{cc}0 & 1 \\0 & 0\end{array}\right); \qquad
Y=\left(\begin{array}{cc}0 & 0 \\1 & 0\end{array}\right).
\end{equation}
We denote by $(h,x,y)\in(\g{sl}(2,\mathbb{R})^*\big)^3$ the dual basis of $(H,X,Y)$.
Thus:
\begin{equation}
\left(\begin{array}{cc}h & x \\ y & -h\end{array}\right)\in\g{sl}(2,\mathbb{R})^*\tens\g{sl}(2,\mathbb{R})
\end{equation}
is the generic point of $\g{sl}(2,\mathbb{R})$.
Let $\NC$ be the nilpotent cone of $\g{sl}(2,\mathbb{R})$. It is the union of three orbits:
\begin{enumerate}
\item $\{0\}$.
\item the  half cone $\NC^+$ with equations $h^2+xy=0;$ $x-y>0$.
\item the  half cone $\NC^-$ with equations $h^2+xy=0;$ $x-y<0$. 
\end{enumerate}
\bigskip

We denote by $\square$ the Casimir operator of $\g{sl}(2,\mathbb{R})$:
\begin{equation}
\square=\frac12 (\de_{H})^2+2\de_{Y}\de_{X}.
\end{equation}
It is an invariant differential operator with constant coefficients on $\g{sl}(2,\mathbb{R})$.

\bigskip

Let $V_{1}=\mathbb {R}^2$ be the standard representation of $\g{sl}(2,\mathbb{R})$.
We denote by $\big(e=(1,0),f=(0,1)\big)$ the standard basis of $\mathbb {R}^2$. The symplectic form $B$ such that $B(e,f)=1$ is $\g{sl}(2,\mathbb{R})$-invariant. For $v\in V_{1}$, we define $\mu_{1}(v)\in\g{sl}(2,\mathbb{R})$ as the unique element such that:

\begin{equation}
\forall Z\in\g{sl}(2,\mathbb{R}),\tr(\mu_{1}(v)Z)=\frac12 B(v,Zv).
\end{equation}
It defines a (moment) map:
\begin{equation}
\mu_{1}:V_{1}\to \g{sl}(2,\mathbb{R}).
\end{equation}
We have
$\mu_{1}(e)=\frac12 X$ and $ \mu_{1}(f)=-\frac12 Y.$
The function $ \mu_{1}$ is a two-fold covering of $\NC^+$ by $V_{1}\setminus\{0\}$.

\bigskip

Let $Z_{0}\in\NC\setminus\{0\}$. Let $\UC$ be a ``small'' neighborhood  of $Z_{0}$.
In this section we determine:
\begin{equation}
\{\phi\in\CC^{-\oo}(\UC,V)^{\g{sl}(2,\mathbb{R})}\,/\, \phi|_{\UC\setminus\NC}=0\}.
\end{equation}
We can assume that $Z_{0}=X\in\NC^+$.

\subsection{Restriction to $X+\mathbb {R}Y$}\label{sec:restr}

We define a map:
\begin{equation}
\begin{split}
\pi:SL(2,\mathbb {R})\times (X+\mathbb {R}Y)&\to \g{sl}(2,\mathbb{R})\\
(g,Z)&\mapsto Ad(g)(Z).
\end{split}
\end{equation}
This map is submersive. Let $I_{2}$ be the identity matrix in $SL(2,\mathbb {R})$. 
Let $\Delta_{X}\subset X+\mathbb {R}Y$ be an open interval containing $X$.
We choose  a connected open subset $\VC\subset SL(2,\mathbb {R})$ such that $I_{2}\in\VC$. We put:
\begin{equation}
\UC=\pi(\VC\mult\Delta_{X}).
\end{equation}
It is an open neighborhood of $X$ in $\g g$.

\begin{lemme}
There is an injective (restriction) map:
\begin{equation}
\begin{split}
\g I_{X}:\CC^{-\oo}(\UC,V)^{\g{sl}(2,\mathbb{R})}&\to \CC^{-\oo}(\Delta_{X},V)\\
\phi&\mapsto\phi_{X}.
\end{split}
\end{equation}
\end{lemme}

\begin{proof}

The map
\begin{equation}
\pi_{\UC}=\pi|_{\VC\mult\Delta_{X}}:\VC\mult\Delta_{X}\to \UC
\end{equation}
is a submersion. Thus if $\phi\in\CC^{-\oo}(\UC,V)$, then $\pi_{\UC}^*(\phi)$ is a well defined generalized function on $\VC\times \D_{X}$ with values in $V$. Moreover,
\begin{equation}
\phi=0\Leftrightarrow \pi_{\UC}^*(\phi)=0.
\end{equation}

Now, we assume that $\phi$ is locally invariant. Then, $\pi_{\UC}^*(\phi)$ is also locally invariant and 
\begin{equation}
\pi_{\UC}^*(\phi)\in\CC^{\oo}(\VC)\what\tens\,  \CC^{-\oo}(\D_{X}).
\end{equation}
(Where $\what\tens$ is a completed tensor product.)
Thus $\pi_{\UC}^*(\phi)$ can be restricted to $\{I_{2}\}\mult \Delta_{X}\subset \VC\mult\Delta_{X}$ (cf. \cite{HC64b}). We identify $\Delta_{X}$ and $\{I_{2}\}\mult \Delta_{X}$. We put:

\begin{equation}
\phi_{X}\defn\pi_{\UC}^*(\phi)|_{ \Delta_{X}}.
\end{equation}
Since $\VC$ is connected and $\phi$ is locally invariant, we have:
\begin{equation}
\pi_{\UC}^{*}(\phi)(g,Z)=\r(g)\phi_{X}(Z).
\end{equation}
Thus 
\begin{equation}
\phi_{X}=0\Leftrightarrow \pi_{\UC}^*(\phi)=0.
\end{equation}
\end{proof}

\bigskip

We have for $Z\in\g{sl}(2,\mathbb{R})$:
\begin{equation}
\LC_{Z}=-h\de_{[Z,H]}-x\de_{[Z,X]}-y\de_{[Z,Y]}.
\end{equation}
In particular:
\begin{eqnarray}
\LC_{H}&=&-2x\de_{X}+2y\de_{Y};\\
\LC_{X}&=&2h\de_{X}-y\de_{H};\\
\LC_{Y}&=&x\de_{H}-2h\de_{Y}.
\end{eqnarray}
If $\VC$ is sufficiently small,  we have $x\not=0$ on $\UC$. We assume that this condition is realized. It follows that on $\UC$ we have:
\begin{equation}\label{eq:radial}
\begin{split}
\de_{X}&=-\frac{1}{2x}\LC_{H}+\frac{y}{x}\de_{Y};\\
\de_{H}&=\frac{1}{x}\LC_{Y}+\frac{2h}{x}\de_{Y}.
\end{split}
\end{equation}

\bigskip

We have $\D_{X}\subset\{X+yY\,/\,y\in\mathbb {R}\}$. We use  the coordinate $y|_{\Delta_{X}}$, still denoted by $y$, on $\Delta_{X}$. Let $\psi\in\CC^{-\oo}(\Delta_{X},V_{n})$. We put $\psi(y)=\psi(X+yY)$.

\begin{lemme}\label{lemme:I_{X}bij}
We have:

\begin{equation}
\g I_{X}\Big(\CC^{-\oo}(\UC,V)^{\g{sl}(2,\mathbb{R})}\Big)=\{\psi\in\CC^{-\oo}(\Delta_{X},V)\,/\, \big(\r(X)+y\r(Y)\big)\psi(y)=0\}.
\end{equation}

Thus:
\begin{equation}
\g I_{X}:\CC^{-\oo}(\UC,V)^{\g{sl}(2,\mathbb{R})}\to\{\psi\in\CC^{-\oo}(\Delta_{X},V)\,/\,
\big(\r(X)+y\r(Y)\big)\psi(y)=0\}
\end{equation}
is an isomorphism

\end{lemme}

\begin{proof}

Since $x|_{\Delta_{X}}=1$ and $h|_{\Delta_{X}}=0$ we have for $\phi\in\CC^{-\oo}(\UC,V)^{\g{sl}(2,\mathbb{R})}$:
\begin{equation}
\begin{split}
\big(\LC_{X}\phi\big)_{X}(y)&=-y\big(\de_{H}\phi\big)_{X}(y);\\
\text { and }\qquad
\big(\LC_{Y}\phi\big)_{X}(y)&=\big(\de_{H}\phi\big)_{X}(y).
\end{split}
\end{equation}
 It follows that we have:
\begin{equation}
\big(\LC_{X}\phi\big)_{X}(y)+y\big(\LC_{Y}\phi\big)_{X}(y)=0.
\end{equation}

\medskip

Let $\psi\in\g I_{X}\Big(\CC^{-\oo}(\UC,V)^{\g{sl}(2,\mathbb{R})}\Big)$. Then,  there is $\phi\in\CC^{-\oo}(\UC,V)^{\g{sl}(2,\mathbb{R})}$ such that $\psi=\phi_{X}$. We have:
\begin{equation}
\begin{split}
\big(\r(X)+y\r(Y)\big)\psi(y)&=\big(\r(X)+y\r(Y)\big)\phi_{X}(y)\\
&=(\r(X)\phi)_{X}(y)+y(\r(Y)\phi)_{X}(y)
+\big(\LC_{X}\phi\big)_{X}(y)+y \big(\LC_{Y}\phi\big)_{X}(y)\\
&=\Big((\r(X)+\LC_{X})\phi\Big)_{X}(y)+y\Big((\r(Y)+\LC_{Y})\phi\Big)_{X}(y)=0.
\end{split}
\end{equation}
\medskip

Let $\psi\in\CC^{-\oo}(\Delta_{X},V) $ such that $\big(\r(X)+y\r(Y)\big)\psi(y)=0.$ We define $\tilde\psi\in\CC^{-\oo}(\VC\times\Delta_{X})$ by the formula:
\begin{equation}
\tilde\psi(g,y)=\r(g)\psi(y).
\end{equation}
Since $\r$ is a smooth function on $SL(2,\mathbb{R})$ with values in $GL(V)$, this is a well defined generalized function on $\VC\mult \Delta_{X}$ with values in $V$. 

Let $(g,Z)\in \VC\times \Delta_{X}$. Let $(g',Z')\in \VC\times \Delta_{X}$ such that $\Ad(g)(Z)=\Ad(g')(Z')$. Then,
$\Ad((g')^{-1}g)Z=Z'$. 
We put $G^{Z}=\{g''\in SL(2,\mathbb {R})\,/\, \Ad(g'')(Z)=Z\}$.  For $g''\in SL(2,\mathbb {R})$, we have
$\Ad(g'')(Z)\in\Delta_{X}\Leftrightarrow  g''\in G^{Z}$. Then, the fiber   of $\pi_{\UC}$ at $(g,Z)$ is included in $\{(g',Z)\,/\, g^{-1}g'\in G^{Z}\}$. Moreover, for $Z'\in\g{sl}(2,\mathbb{R})$, $[Z,Z']=0\Leftrightarrow Z'\in\mathbb {R}Z$. Thus, since $\VC$ is connected, the condition $\big(\r(X)+y\r(Y)\big)\psi(y)=0$ on $\Delta_{X}$ ensures that $\tilde\psi$ is constant along the fibers of $\pi_{\UC}$. Thus there is a well defined generalized function $\overline\psi$ on $\UC$ such that:
\begin{equation}
\pi_{\UC}^*(\overline \psi)=\tilde\psi.
\end{equation}
It follows from the construction that 
$\big(\overline\psi\big)_{X}=\psi.$

\end{proof}

\bigskip

The hypothesis $\phi|_{\UC\setminus\NC}=0$ means that $\phi_{X}$ is supported in $\{X\}\subset\Delta_{X}$.

\subsection{Radial part of $\square$}
In the neighborhood $\UC$ of $X$ defined in section \ref{sec:restr}:
\begin{equation}
\begin{split}
\square&=\frac{1}{2}(\de_{H})^2+2\de_{Y}\de_{X}\\
&=\frac{1}{2}\big(\frac{1}{x}\LC_{Y}+\frac{2h}{x}\de_{Y}\big)^2+2\de_{Y}\big(\frac{-1}{2x}\LC_{H}+\frac{y}{x}\de_{Y}\big).
\end{split}
\end{equation}
We define the radial part of $\square$ as the differential operator $\square_{X}$ on $\CC^{-\oo}(\Delta_{X},V)$:
\begin{equation}
\begin{split}
\square_{X}&=\Big(3+\r(H)+2y\frac{\de}{\de y}\Big)\frac{\de}{\de y}+\frac{1}{2}\r(Y)^2.
\end{split}
\end{equation}
This definition is justified by the following lemma:

\begin{lemme}
Let $\phi\in\CC^{-\oo}(\UC,V)^{\g{sl}(2,\mathbb{R})}$, then we have:

\begin{equation}
(\square\phi)_{X}=\square_{X}\phi_{X}.
\end{equation}

\end{lemme}

\begin{proof}
Since $x|_{\Delta_{X}}=1$ and $h|_{\Delta_{X}}=0$, we have:
\begin{equation}
\begin{split}
(\square\phi)_{X}&=\frac12\Big(\big(\LC_{Y}^2+2\LC_{Y}\frac{h}{x}\frac{\de}{\de y}\big)\phi\Big)_{X}
+2\Big(\frac{-1}{2}\big(\frac{\de}{\de y}\LC_{H}\phi\big)_{X}+\big(\frac{\de}{\de y}y\frac{\de}{\de y}\phi\big)_{X}\Big)\\
&=\frac12\Big(\big(\r(Y)^2+2(x\de_{H}-2h\de_{Y})\frac{h}{x}\frac{\de}{\de y}\big)\phi\Big)_{X}\\
&\hskip 2cm+2\Big(\frac{-1}{2}\big(-\r(H)\frac{\de}{\de y}\phi_{X}\big)+\frac{\de}{\de y}\phi_{X}+y\big(\frac{\de}{\de y}\big)^2\phi_{X}\Big)\\
&=\frac12\Big(\r(Y)^2+2\frac{\de}{\de y}\Big)\phi_{X}
+\Big(\r(H)+2+2y\frac{\de}{\de y}\Big)\frac{\de}{\de y}\phi_{X}\\
&=\Big(3+\r(H)+2y\frac{\de}{\de y}\Big)\frac{\de}{\de y}\phi_{X}+\frac{1}{2}\r(Y)^2\phi_{X}=\square_{X}\phi_{X}.
\end{split}
\end{equation}

\end{proof}

\subsection{The Dirac function $\d_{\NC^{+}}$ (resp. $\d_{\NC^{-}}$)}
\label{sec:Dirac}
 
 Let $dZ=dx\,dy\,dh$ be the Lebesgue measure on $\g{sl}(2,\mathbb{R})$. 
 Let $(e^*,f^*)\in (V_{1}^*)^2$ be the dual basis of  $(e,f)$.
The Lebesgue measure $dv=-2 de^*\,df^*$ on $V_{1}$ is $\g{sl}(2,\mathbb{R})$-invariant. We define an invariant generalized function  $\d_{\NC^+}$ (resp. $\d_{\NC^{-}}$)  on $\g{sl}(2,\mathbb{R})$ and supported in $\NC^+\cup\{0\}$ (resp. $\NC^-\cup\{0\}$) by:
\begin{equation}
 \begin{split}
\forall g\in\CC^{\oo}_{c}(\g{sl}(2,\mathbb{R})),\quad
\int_{\g{sl}(2,\mathbb{R})}\d_{\NC^+}(Z)g(Z)dZ&\mathop{=}\limits^{\text{def}}\int_{V_{1}}g\circ\mu_{1}(v)dv\\
\Big(\text{resp.}
\forall g\in\CC^{\oo}_{c}(\g{sl}(2,\mathbb{R})),\quad
\int_{\g{sl}(2,\mathbb{R})}\d_{\NC^-}(Z)g(Z)dZ&\mathop{=}\limits^{\text{def}}\int_{V_{1}}g\circ(-\mu_{1})(v)dv
\Big).
\end{split}
\end{equation}
We put:
\begin{equation}
\d_{X}=(\d_{\NC^+})_{X}\in\CC^{-\oo}(\Delta_{X}).
\end{equation}
 We still denote by $dy$ the Lebesgue measure on $\Delta_{X}$. It is invariant. 
Let $g\in\CC^{\oo}_{c}(\Delta_{X})$. Then we have:
\begin{equation}
\int_{\D_{X}}\d_{X}(y)g(y)dy=g(0).
\end{equation}

\subsection{Irreducible representations }\label{sec:IrrRep}

If $V=V^{1}\oplus\dots\oplus V^{n}$ where $V^{i}$ is  an irreducible representation of $\g{sl}(2,\mathbb{R})$, then we have:
\begin{equation}
\CC^{-\oo}(\UC,V)=\sdir_{i=1}^n\CC^{-\oo}(\UC,V^{i}),
\end{equation}
every subspace being stable for $\g{sl}(2,\mathbb{R})$. Thus we can assume from now on that the representation of  $\g{sl}(2,\mathbb{R})$ in $V$ is irreducible.

\medskip

We fix the Cartan subalgebra $\g h=\mathbb {R}H$ and the positive root $2h$ (we still denote by $h$ its restriction to ${\g h}$).
Let $n\in\mathbb {N}$. We denote by $V_{n}$ the irreducible   representation of $\g{sl}(2,\mathbb{R})$ with highest weight $nh$.
We have $\dim(V_{n})=n+1$.
We decompose $V_{n}$ under the action of $\mathbb{R}H$.
We fix $v_{0}\in V_{n}\setminus \{0\}$ a vector of weight $-nh$:
\begin{equation}
\r(H)v_{0}=-nv_{0}.
\end{equation}
We put for $0\le i\le n$: $v_{i}=\r(X)^{i}v_{0}$. We have $\r(X)v_{n}=0$ and
$\r(H)v_{i}=(-n+2i)v_{i}.$
On the other hand, $\r(Y)v_{0}=0$ and for $1\le i\le n$:
$\r(Y)v_{i}=(n-i+1)iv_{i-1}.$

\subsection{A basic function on $\NC^+$}
\label{sec:s_n}
We construct a function  $s_{n}:\UC\cap\NC^+\to V_{n}$ which is the basic tool to generate all the generalized functions we are looking for.

\subsubsection{Case $n$ even} In this case $V_{n}$ is isomorphic to the irreducible component of $S^{\frac{n}{2}}(\g{sl}(2,\mathbb{R}))$  (under adjoint action of $\g{sl}(2,\mathbb{R})$) generated by $X^{\frac{n}{2}}$.
From now on we will identify $V_{n}$ with this component.
We denote by $s_{n}:\NC\to V_{n}$ the invariant map defined by:
\begin{equation}
s_{n}(Z)=Z^{\frac{n}{2}}.
\end{equation}

\subsubsection{Case $n=1$} 

We recall that $\mu_{1}:V_{1}\setminus\{0\}\to \NC^+$ is a two-fold covering with $\mu_{1}(e)=\frac12 X$.
If $\UC$ is  a sufficiently small connected neighborhood of $X$, there exists a unique continuous section $s_{1}$ of $\mu_{1}$ in $\UC\cap\NC^+$ such that $s_{1}(\frac12 X)=e$. We have $s_{1}:\UC\cap\NC^+\to V_{1}$. It satisfies:
\begin{equation}
\forall Z\in\UC\cap\NC^+,\ \mu_{1}(s_{1}(Z))=Z.
\end{equation}

\subsubsection{Case $n$ odd} 
More generally, when $n$ is odd, $V_{n}$ is isomorphic to the irreducible component of $V_{1}\tens S^{\frac{n-1}{2}}(\g{sl}(2,\mathbb{R}))$ generated by 
$e\tens X^{\frac{n-1}{2}}$.
From now on we will identify $V_{n}$ with this component.
Let $\UC$ be the above neighborhood  of $X$. We define a function $s_{n}:\UC\cap\NC^+\to V_{n}$ by:
\begin{equation}
\forall Z\in\UC\cap\NC^+,\ s_{n}(Z)=s_{1}(Z)\tens Z^{\frac{n-1}{2}}\in V_{n}.
\end{equation}

\subsection{Basic theorem}
Let $\UC$ be an open subset of $\g{sl}(2,\mathbb{R})$. We put:
\begin{equation}
\SC_{n}(\UC)=\big\{\phi\in\CC^{-\oo}(\UC,V_{n})^{\g{sl}(2,\mathbb{R})}\,/\, \phi|_{\UC\setminus\NC}=0\}.
\end{equation} 

\begin{theo}\label{th:local}
Let $n\in\mathbb{N}$.  Let $\UC$ be an open connected neighborhood  of $X$ such that the  function $s_{n}$ is well defined on $\UC\cap\NC$ (cf. section \ref{sec:s_n}) and $\g I_{X}$ is bijective (cf. section \ref{sec:restr}). 
Then:
\begin{enumerate}
\item When $n$ is even, $\SC_{n}(\UC)$ is an infinite dimensional vector space  with basis:
\begin{equation}
\big(\square^k(s_{n}\d_{\NC^+})\big)_{k\in\mathbb{N}}.
\end{equation}

\item When $n$ is odd, $\dim(\SC_{n}(\UC))=\frac{n+1}{2}$ and a basis is given by:
\begin{equation}
\big(\square^k(s_{n}\d_{\NC^+})\big)_{0\le k\le \frac{n-1}{2}}.
\end{equation}

\end{enumerate}

\end{theo}

{\bf Remark:}
Since $\d_{\NC^+}(Z)dZ$ is a measure on $\g{sl}(2,\mathbb{R})$ with support $\NC^+\cup\{0\}$ and $s_{n}$ is a smooth function on $\UC\cap\NC$ with values in $V_{n}$, $s_{n}\d_{\NC^+}$ is a well defined generalized function on $\UC$ with values in $V_{n}$.

\begin{proof} 
Thanks to the isomorphism $\g I_{X}$ we have to  determine the space:
\begin{equation}
\{\psi\in\CC^{-\oo}(\Delta_{X},V_{n})\,/\, \psi|_{\Delta_{X}\setminus\{0\}}=0\text{ and }\big(\r(X)+y\r(Y)\big)\psi(y)=0\}.
\end{equation}
Let $\psi\in\CC^{-\oo}(\Delta_{X},V_{n})$. We write:
\begin{equation}
\psi(y)=\som_{i=0}^n \psi_{i}(y)v_{i},
\end{equation}
where $\psi_{i}\in\CC^{-\oo}(\Delta_{X})$ and $(v_{i})_{0\le i\le n}$ is the basis defined in section \ref{sec:IrrRep}. 
We put:
\begin{equation}
\d^{k}(y)=\Big(\frac{\de}{\de y}\Big)^k\d_{X}(y).
\end{equation}
Since $\psi$  is supported in $\NC$ and $\Delta_{X}\cap \NC=\{X\}$, there exists $a_{i,k}\in\mathbb{R}$, all equal to zero but for finite number, such that:
\begin{equation}
\psi_{i}(y)=\som_{k\in\mathbb{N}}a_{i,k}\d^{k}(y).
\end{equation}

For $n=0$,  we have $\r=0$ and the condition $\big(\r(X)+y\r(Y)\big)\psi(y)=0$ is automatically satisfied. 

For $n\ge 1$, we put $\a_{i}=(n-i+1)i$. We have $y\d^{0}(y)=0$ and for $k\ge 1$, $y\d^{k}(y)=-k\d^{k-1}(y)$. Thus:
\begin{equation}
\som_{0\le i\le n-1,\,k\in\mathbb{N}}a_{i,k}\d^{k}(y)v_{i+1}
-\som_{1\le i\le n,\,k\ge 1}\a_{i}a_{i,k}k\d^{k-1}(y)v_{i-1}=0.
\end{equation}
It follows:
\begin{equation}
\begin{cases}
a_{n-1,k}=0&\text{for }k\ge 0;\\
a_{1,k}=0&\text{for }k\ge 1;\\
a_{i-1,k}=(k+1)(i+1)(n-i)a_{i+1,k+1}& \text{for }n\ge 2,\ 1\le i\le n-1\text{ and }k\ge 0.
\end{cases}
\end{equation}
It follows in particular
\begin{enumerate}
\item \label{rel:i} from the first and the last relations that $\forall i,k\ge 0\,$ with $2i+1\le n$: $a_{n-(2i+1),k}=0$;
\item \label{rel:ii}
from the last relation that $\forall i\ge 0$ with $2i\le n$, $(a_{n-2i,k})_{k\ge 0}$ is completely determined by $(a_{n,k})_{k\ge 0}$.
\end{enumerate}

\medskip

We distinguish between the two cases according to the parity of $n$.

\medskip

\emph{$n$ even:}
In this case, for $n\ge 2$, the second relation follows from (\ref{rel:i}). Hence the map:
\begin{equation}\label{eq:IsoEven}
\begin{split}
\{\psi\in\CC^{-\oo}(\Delta_{X},V_{n})\,/\, \psi|_{\Delta_{X}\setminus\{0\}}=0\text{ and }\big(\r(X)+y\r(Y)\big)\psi(y)=0\} &\to \mathbb{R}^{\mathbb {N}}\\
\psi(y)=\som_{0\le i\le n,\, k\in\mathbb{N}}a_{i,k}\d^{k}(y)v_{i}&\mapsto (a_{n,k})_{k\in\mathbb{N}}
\end{split}
\end{equation}
is bijective. This is also true for $n=0$.

\medskip

\emph{$n$ odd:}
It follows from the two last relations that for $k\ge i\ge 1$
$a_{2i-1,k}=0.$
In particular, the map:
\begin{equation}\label{eq:IsoOdd}
\begin{split}
\{\psi\in\CC^{-\oo}(\Delta_{X},V_{n})\,/\,\psi|_{\Delta_{X}\setminus\{0\}}=0\text{ and }\big(\r(X)+y\r(Y)\big)\psi(y)=0\} &\to \mathbb{R}^{\frac{n+1}{2}}\\
\psi(y)=\som_{0\le i\le n,\, k\in\mathbb{N}}a_{i,k}\d^{k}(y)v_{i}&\mapsto (a_{n,0},\dots,a_{n,\frac{n-1}{2}})
\end{split}
\end{equation}
is bijective.

\bigskip

This proves the first part of the theorem on the dimension of $\SC_{n}(\UC)$. It remains to prove that the  functions $\square^k(s_{n}\d_{\NC})$ form a basis of $\SC_{n}(\UC)$.
 We have for $\psi(y)=\som_{i=0}^n \som_{k\in\mathbb {N}}a_{i,k}
\d^{k}(y)v_{i}\in\CC^{-\oo}(\Delta_{X},V_{n})$ such that $\r(X+yY)\psi(y)=0$:
 \begin{equation}
\begin{split}
\square_{X}\psi(y)
&=
\big(3+
\r(H)+2y\de_{Y}\big)\som_{k\in\mathbb {N}}a_{n,k}\d^{k+1}(y)v_{n}+\som_{i=0}^{n-1} \,\dots\,v_{i}\\
&=
\som_{k\in\mathbb {N}}(n-2k-1)a_{n,k} \d^{k+1}(y)v_{n}
+\som_{i=0}^{n-1} \,\dots\,v_{i}
\end{split}
\end{equation}
where $\dots$ are elements of $\CC^{-\oo}(\Delta_{X})$.

\medskip

\emph{$n$ even:} Since $v_{n}=X^{\frac{n}{2}}$, we have
$(s_{n}\d_{\NC})_{X}(y)=\d_{X}(y)X^{\frac{n}{2}}.$
By induction on $k$, it follows:
\begin{equation}
\begin{split}
\big(\square^k(s_{n}\d_{\NC})\big)_{X}(y)&=(n-2k+1)\dots(n-1)\d^{k}(y)X^{\frac{n}{2}}\\
&\hskip 2cm +\text{ terms with } X^{\frac{n}{2}-i} \text{ for } i\ge 1.
\end{split}
\end{equation}
Since $n$ is even $n-2k+1\not=0$. The result follows.

\emph{$n$ odd:} Since $v_{n}=e\tens X^{\frac{n-1}{2}}$, we have
$(s_{n}\d_{\NC})_{X}(y)=\d_{X}(y)(e\tens X^{\frac{n-1}{2}}).$
By induction on $k$, it follows:
\begin{equation}
\begin{split}
\big(\square^k(s_{n}\d_{\NC})\big)_{X}(y)&=(n-2k+1)\dots(n-1)\d^{k}(y)(e\tens X^{\frac{n-1}{2}})\\
&\hskip 2cm +\text{ terms with } e\tens X^{\frac{n-1}{2}-i} \text{ for } i\ge 1.
\end{split}
\end{equation}
In this case for $k=\frac{n+1}{2}$, $n-2k+1=0$. Thus, since $\square^k (s_{n}\d_{\NC})$ is invariant, it follows from the isomorphism (\ref{eq:IsoOdd}) that for $k\ge \frac{n+1}{2}$:
$\square^k (s_{n}\d_{\NC})=0.$
The result follows.

\end{proof}

\bigskip

\subsection{Global version}

Let $\UC$ be an open subset of $\g{sl}(2,\mathbb{R})$.
We put:
\begin{eqnarray}
\SC_{n}^{0}(\UC)&=&\{\phi\in\CC^{-\oo}(\UC,V_{n})^{\g{sl}(2,\mathbb{R})}\,/\, \phi|_{\UC\setminus\{0\}}=0\};\\
\SC_{n}^{\pm}(\UC)&=&\{\phi\in\CC^{-\oo}(\UC,V_{n})^{\g{sl}(2,\mathbb{R})}\,/\, \phi|_{\UC\setminus(\NC^{\pm}\cup\{0\})}=0\}.
\end{eqnarray}

\bigskip

\begin{theo}\label{th:global}
Let $\UC$ be an $SL(2,\mathbb {R})$-invariant open subset of $\g{sl}(2,\mathbb{R})$. Then we have:

\begin{enumerate}
\item 
\begin{equation}
\begin{cases}
\SC_{n}^{0}(\UC)=\{0\}&\text {if }0\not\in\UC;\\
\SC_{n}^{0}(\UC)\simeq (V_{n}\tens S(\g{sl}(2,\mathbb{R})))^{\g{sl}(2,\mathbb{R})}
&\text{if }0\in\UC.
\end{cases}
\end{equation}

\item When $n$ is even, we have:
\begin{eqnarray}
\SC_{n}(\UC)&=&\SC_{n}^{0}(\UC)\oplus\Vect\{\square^k(s_{n}\d_{\NC^{+}})|_{\UC}/k\in\mathbb{N}\}
\oplus\Vect\{\square^k( s_{n}\d_{\NC^{-}})|_{\UC}/k\in\mathbb{N}\}
\\
\SC_{n}^{\pm}(\UC)&=&\SC_{n}^{0}(\UC)\oplus\Vect\{\square^k( s_{n}\d_{\NC^{\pm}})|_{\UC}/k\in\mathbb{N}\}
\end{eqnarray}

\item When $n$ is odd:
\begin{equation}
\SC_{n}(\UC)=\SC_{n}^{\pm}(\UC)=\SC_{n}^{0}(\UC)
\end{equation}

\end{enumerate}

\end{theo}

\begin{proof}

(i) It follows  from Theorem \ref{th:Supp0}.

(ii)
When $n$ is even, the function $\d_{\NC}^{\pm}$ is defined on $\g{sl}(2,\mathbb{R})$, the function $s_{n}$ is  defined on $\NC$ and the product $s_{n}\d_{\NC^{\pm}}$ is well defined (cf. Remark of Theorem \ref{th:local}). Then the result  follows from Theorem \ref{th:local}.

(iii) 
 Let $n$ be odd. We assume that $\UC\cap \NC\not=\emptyset$. Since $\UC$ is $SL(2,\mathbb {R})$-invariant, we have $\NC^+\subset\UC$ or $\NC^{-}\subset \UC$. We assume that $\NC^+\subset\UC$ (the case $\UC\subset \NC^-$ is similar).

Let $\phi\in\CC^{-\oo}(\UC,V)^{\g{sl}(2,\mathbb{R})}$. Let $\UC_{0}\subset \UC$ be a suitable neighborhood  of $X$ where $s_{1}$ (and thus $s_{n}$) is defined (cf. section \ref{sec:s_n}). There exists $(a_{0},\dots ,a_{\frac{n-1}2})\in\mathbb {R}^{\frac{n+1}2}$ such that  on $\UC_{0}$ (cf. Theorem \ref{th:local}):
\begin{equation}
\phi(Z)=\som_{k=0}^{\frac{n+1}2}a_{k}\square^k\big(s_{n}(Z)\d_{\NC^+}(Z)\big)
=\som_{k=0}^{\frac{n+1}2}a_{k}\square^k\big((s_{1}(Z)\tens Z^{\frac{n-1}2})\d_{\NC^+}(Z)\big).
\end{equation}

Since $\mu_{1}:V_{1}\setminus \{0\}\to \NC^+$ is a non trivial two-fold covering, there is not any  continuous section. In other words there is not any continuous $SL(2,\mathbb {R})$-invariant map $s:\NC^+\to V_{1}$ such that for any  $Z\in\UC_{0}$, $s(Z)=s_{1}(Z)$. Thus $a_{0}=\dots=a_{\frac{n-1}2}=0$. The result follows.

\end{proof}

\section{Invariant solutions of differential equations}

\subsection{Introduction}
Let $\CC_{V}=(\End(V)\tens S(\g{sl}(2,\mathbb{R})))^{\g{sl}(2,\mathbb{R})}$ be the algebra of $\End(V)$-valued invariant differential operators with constant coefficients  on $\g g$. It is the \emph{classical family algebra} in the terminology of Kirillov (cf. \cite{Kir00}). When $V=V_{n}$ is the $(n+1)$-dimensional irreducible representation of $\g{sl}(2,\mathbb{R})$, we put $\CC_{n}=\CC_{V_{n}}$.

Let $\UC\subset\g{sl}(2,\mathbb{R})$ be an open subset. It is a natural and interesting problem to determine  the generalized functions $\phi\in\CC^{-\oo}(\UC,V)^{\g{sl}(2,\mathbb{R})}$ such that $\CC_{V}\phi$ is finite dimensional.

\bigskip

We recall that $S(\g{sl}(2,\mathbb{R}))^{\g{sl}(2,\mathbb{R})}=\mathbb {R}[\square]$. It is a subalgebra of $\CC_{V}$. An other subalgebra of $\CC_{V}$ is $\End(V)^{\g{sl}(2,\mathbb{R})}$. 
When $V=V_{n}$, we put:
\begin{equation}
M_{n}=\r_{n}(X)Y+\r_{n}(Y)X+\frac12 \r_{n}(H)H\in \CC_{n}
\end{equation}
 According N. Rozhkovskaya (cf. \cite{Roz03}), $\CC_{n}$ is a free $S(\g{sl}(2,\mathbb{R}))^{\g{sl}(2,\mathbb{R})}$-module with basis $\BC_{n}=\big(1,M_{n},\dots, (M_{n})^n\big).$

\begin{lemme}
Let $\phi\in\CC^{-\oo}(\UC,V)^{\g{sl}(2,\mathbb{R})}$. Then we have:
\begin{equation}
\dim_{\mathbb {R}}\big(\CC_{V}\phi\big)<\oo\Leftrightarrow
\dim_{\mathbb {R}}\big(\mathbb {R}[\square]\phi\big)<\oo
\end{equation}

\end{lemme}

\begin{proof}
We argue as in \cite{Roz03}.
Let $H$ be the set of harmonic polynomials in $S(\g{sl}(2,\mathbb{R}))$. Then, $S(\g{sl}(2,\mathbb{R}))=\mathbb {R}[\square]\tens H$  (cf. \cite{Kos}), and:
\begin{equation}
\CC_{V}=\mathbb {R}[\square]\tens\big(H\tens \End(V)\big)^{\g{sl}(2,\mathbb{R})}.
\end{equation}
Since $\dim_{\mathbb {R}}\big(H\tens \End(V)\big)^{\g{sl}(2,\mathbb{R})}<\oo$, the result follows.

\end{proof}

{\bf Remark:}
Since $\mathbb {R}[\square]\subset\mathbb {R}[\square]\tens\End(V)^{\g{sl}(2,\mathbb{R})}\subset\CC_{V}$, the condition $\dim(\CC_{V}\phi)<\oo$ is also equivalent to the existence of $r\in\mathbb {N}$ and $(A_{0},\dots,A_{r-1})\in\big(\End(V)^{\g{sl}(2,\mathbb{R})}\big)^{r}$ such that:
\begin{equation}\label{eq:generique0}
\big(\square^r+A_{r-1}\square^{r-1}+\dots A_{1}\square+A_{0}\big)\phi=0.
\end{equation}
Useful examples of (\ref{eq:generique0}) are $(\square-\l)^k\phi=0$ for $\l\in\mathbb {C}$ and  generalized functions with values in a complex representation. We give such an example below.

\begin{defi}
Let $\phi\in\CC^{-\oo}(\UC,V)^{\g{sl}(2,\mathbb{R})}$. We say that $\phi$ is $\square$-finite if $\dim_{\mathbb {R}}(\mathbb {R}[\square]\phi)<\oo$.
\end{defi}
In other words, $\phi$ is $\square$-finite if there exists $r\in\mathbb {N}$ and $(a_{0},\dots,a_{r-1})\in\mathbb {R}^{r}$ such that
\begin{equation}\label{eq:generique}
\big(\square^r+a_{r-1}\square^{r-1}+\dots a_{1}\square+a_{0}\big)\phi=0.
\end{equation}

\bigskip

{\bf Example:} (This was our original motivation to study this problem.) Let $\g g=\g g_{\0}\oplus\g g_{\1}$ be a Lie superalgebra. We define the generalized functions on  $\g g$ as the generalized functions on $\g g_{\0}$ with values in the exterior algebra $\L(\g g_{\1}^*)$ of $\g g_{1}^*$:
\begin{equation}
\CC^{-\oo}(\g g)\defn \CC^{-\oo}(\g g_{\0})\tens \L(\g g_{\1}^*)=\CC^{-\oo}(\g g_{\0},\L(\g g_{1}^*)).
\end{equation}
We assume that $\g g$ has a non degenerate invariant symmetric even bilinear form $B$. Let $\W\in S^{2}(\g g)$ be the Casimir operator associated with $B$. We have $\W=\W_{\0}+\W_{\1}$ with $\W_{\0}\in S^2(\g g_{\0})$ and $\W_{\1}\in\L^2(\g g_{\1})$. We consider $\W_{\1}$ as an element of $\End(\L(\g g_{\1}^*))$ acting by interior product. When they can be evaluated (cf. for example \cite[Chapitre III.5]{Lav98}), the Fourier transforms of the coadjoint orbits in $\g g^*$ are invariant generalized functions $\phi$ on $\g g$ subject to equations of the form $(\W-\l)\phi=0$ with $\l\in\mathbb {C}$. It can be written $\big(\W_{\0}+(\W_{\1}-\l)\big)\phi=0$ (for $\g g_{\0}=\g{sl}(2,\mathbb{R})$ it is of the form (\ref{eq:generique0}) with $\W_{\0}=\square$ and $A_{\0}=\W_{\1}-\l$).
We have:
\begin{equation}
(\W_{\0}-\l)^{k}=\som_{i=0}^{k}\Big(\begin{matrix}k \\ i\end{matrix}\Big)(\W-\l)^{i}(-\W_{\1})^{k-i}.
\end{equation}
For $k> \frac{\dim(\g g_{\1})}{2}$, we have $\W_{\1}^{k}=0$. It follows that for $k>1+ \frac{\dim(\g g_{\1})}{2}$ we have:
\begin{equation}
(\W_{\0}-\l)^{k}\phi=0.
\end{equation}
this equation is of the form of (\ref{eq:generique0}).

\subsection{Generalized functions with support $\{0\}$}

We immediately obtain  from Theorem \ref{th:Supp0}

\begin{theo}\label{th:InvSupp0}  Let $V$ be a representation of $\g{sl}(2,\mathbb{R})$. Let $\phi\in\CC^{-\oo}(\g{sl}(2,\mathbb{R}),V)^{\g{sl}(2,\mathbb{R})}$ such that $\phi|_{\g{sl}(2,\mathbb{R})\setminus\{0\}}=0$ and $\phi$ is $\square$-finite.
Then, we have  $\phi=0$.

\end{theo}

\subsection{Support in the nilpotent cone: local version}

\begin{theo}\label{th:InvLocal}
Let $n\in\mathbb N$. Let $V_{n}$ be the irreducible $n+1$-dimensional representation of $\g{sl}(2,\mathbb{R})$. Let $W$ be a finite dimensional vector space with trivial action of $\g{sl}(2,\mathbb{R})$.   Let $\UC$ be an open connected neighborhood  of $X$ such that the  function $s_{n}$ is well defined on $\UC\cap\NC$ (cf. section \ref{sec:s_n}) and $\g I_{X}$ is bijective (cf. section \ref{sec:restr}).  Let $\phi\in\CC^{-\oo}(\UC,W\tens V_{n})^{\g{sl}(2,\mathbb{R})}$ such that $\phi|_{\UC\setminus \NC}=0$. 
Let $r\in\mathbb {N}$ and $(a_{0},\dots,a_{r-1})\in\mathbb {R}^{r}$
such that:
$\Big(\square^r+\som_{k=0}^{r-1} a_{k}\square^k\Big)\phi=0.$

Then, we have $\phi=0$ when at least one of the following conditions is satisfied:
\begin{enumerate}
\item  $n$ is even;

\item  $n$ is odd and  $a_{0}\not=0$.
\end{enumerate}

\end{theo}

\begin{proof}

Let $\phi\in\CC^{-\oo}(\UC,W\tens V_{n})^{\g{sl}(2,\mathbb{R})}$ such that $\phi|_{\UC\setminus \NC}=0$. From Theorem \ref{th:local} we obtain that there exist  $p\in\mathbb N$, with $p=\frac{n-1}{2}$ if $n$ is odd and $(w_{0},\dots,w_{p})\in W^{p+1}$, such that:

\begin{equation}
\phi=\som_{i=0}^pw_{i}\tens\square^{i}(s_{n}\d_{\NC^+}).
\end{equation}
Then:
\begin{enumerate}
\item When $n$ is even, for $0\le j\le p+r$, we have
$\som_{k+i=j}a_{k}w_{i}=0.$

\item When $n$ is odd, for $0\le j\le \frac{n-1}{2}$, we have
$\som_{k+i=j}a_{k}w_{i}=0.$

\end{enumerate}
The result follows.

\end{proof}

{\bf Remark:} When $n$ is odd, in contrast  with the classical case ($V=V_{0}$ is the trivial representation) there exist (in a neighborhood of $X$) non trivial locally invariant solutions of the equation $\square^k\phi=0$ supported in the nilpotent cone! For example, if $k\ge\frac{n+1}2$ the functions $\phi=\square^{i}(s_{n}\d_{\NC^+})$ for $0\le i\le \frac{n-1}2$ are not trivial, supported in the nilpotent cone and satisfy the equation $\square^{k}\phi=0$.

When we consider the equation  $(\square-\l)^k\phi=0$ for $\l\in\mathbb {C}\setminus\{0\}$, then the trivial solution is again the only one supported in the nilpotent cone.

\subsection{Support in the nilpotent cone: global version}

\begin{theo}\label{th:InvGlobal}
 Let $V$ be a real finite dimensional representation of $\g{sl}(2,\mathbb{R})$. Let $\UC$ be an $SL(2,\mathbb {R})$-invariant open subset of $\g{sl}(2,\mathbb{R})$. 
Let $\phi\in\CC^{-\oo}(\UC,V)^{\g{sl}(2,\mathbb{R})}$ such that $\phi|_{\UC\setminus\NC}=0$ and $\phi$ is $\square$-finite.
Then we have $\phi=0$.
\end{theo}

\begin{proof}
It is enough to prove the theorem for $V$ irreducible. 
Then, the result follows from Theorem \ref{th:global}, Theorem \ref{th:InvLocal} and Theorem \ref{th:InvSupp0}.
\end{proof}

\section{General invariant generalized functions}

\subsection{Main theorem}

\begin{theo}\label{th:Principal}
 Let $V$ be a real finite dimensional representation of $\g{sl}(2,\mathbb{R})$. Let $\UC$ be an $SL(2,\mathbb {R})$-invariant open subset of $\g{sl}(2,\mathbb{R})$. 
Let $\phi\in\CC^{-\oo}(\UC,V)^{\g{sl}(2,\mathbb{R})}$ 
such that  $\phi$ is $\square$-finite.
Then $\phi$ is determined  by $\phi|_{\UC\setminus \NC}$ and  $\phi|_{\UC\setminus \NC}$  is an analytic function. 
\end{theo}

\begin{proof}
The fact that $\phi$ is determined by $\phi|_{\UC\setminus \NC}$ follows from Theorem  \ref{th:InvGlobal}.
The fact that $\phi|_{\UC\setminus \NC}$ is analytic can be proved exactly as in \cite{HC65a}.

\end{proof}

{\bf Remark:}
In general $\phi$ will not be locally $L^1$.
Indeed, let $\phi_{0}\in\CC^{-\oo}(\g{sl}(2,\mathbb{R}))^{\g{sl}(2,\mathbb{R})}$ a non zero $\square$-finite generalized function. Then $\phi_{0}$ is locally $L^{1}$, but for $k\in\mathbb {N}^*$:
\begin{equation}
M_{n}^{k}\phi_{0}\in\CC^{-\oo}\big(\g{sl}(2,\mathbb{R}),\End(V_{n})\big)^{\g{sl}(2,\mathbb{R})}
\end{equation}
is usually not locally $L^1$.

\bigskip

\subsection{Application to the Superpfaffian} 
Let us consider the Lie superalgebra $\g g=\g{spo}(2,2n)$. Its even part is $\g g_{\0}=\g{sl}(2,\mathbb{R})\oplus\g{so}(2n,\mathbb {R})$. Its odd part is $\g g_{\1}= V_{1}\tens W$ where $W$ is the standard $2n$-dimensional representation of $\g{so}(2n,\mathbb {R})$.  

In \cite{Lav04} we constructed a particular   invariant generalized function $\Spf$ on $\g{spo}(2,2n)$ called Superpfaffian. It generalizes the Pfaffian on $\g{so}(2n,\mathbb {R})$ and the inverse square root of the determinant on $\g{sl}(2,\mathbb{R})$. As it is a polynomial of degree $n$ on $\g{so}(2n,\mathbb {R})$, we may consider that we have:
\begin{equation}
\Spf\in\CC^{-\oo}\Big(\g{sl}(2,\mathbb{R})\,,\,\sdir_{k=0}^n S^{k}(\g{so}\big(2n,\mathbb {R})^*\big)\tens \L\big(\g g_{1}^*\big)\Big)^{\g{sl}(2,\mathbb{R})}.
\end{equation}

Let $\W$ (resp. $\square$, $\W'_{\0}$, $\W_{\1}$) be the Casimir operator on $\g{spo}(2,2n)$ (resp. on $\g{sl}(2,\mathbb{R})$, $\g{so}(2n,\mathbb {R})$, $\g g_{1}$). Then $\W=\square+\W'_{\0}+\W_{\1}$ and
\begin{equation}
\W'_{\0}+\W_{\1}\in\End\Big(\sdir_{k=0}^n S^{k}\big(\g{so}(2n,\mathbb {R})^*\big)\tens \L\big(\g g_{1}^*\big)\Big)^{\g{sl}(2,\mathbb{R})}
\end{equation}
is a nilpotent endomorphism.
 The superpfaffian satisfies:
\begin{equation}
\big(\square+(\W'_{\0}+\W_{\1})\big)\Spf=\W\Spf=0.
\end{equation}
The function $\Spf$ is analytic on $\g{sl}(2,\mathbb{R})\setminus\NC$ and  in \cite{Lav04} an explicit formula is given for $\Spf(X)\in\sdir_{k=0}^n S^{k}\big(\g{so}(2n,\mathbb {R})^*\big)\tens \L\big(\g g_{1}^*\big)$ with $X\in\g{sl}(2,\mathbb{R})\setminus\NC$. However,
since  $\Spf$ is not locally $L^1$ (cf. \cite{Lav04}), it is not clear whether   $\Spf$ is determined by  its restriction to $\g{sl}(2,\mathbb{R})\setminus\NC$ or not. In \cite{Lav04} we proved that $\Spf$ is characterized, as an invariant generalized function on $\g{sl}(2,\mathbb{R})$,  by its restriction to $\g{sl}(2,\mathbb{R})\setminus \NC$ and its wave front set. 

From the preceding results we obtain this new characterization of $\Spf$:

\begin{theo}
Let $\phi\in\CC^{-\oo}\Big(\g{sl}(2,\mathbb{R})\, ,\, \sdir_{k=0}^n S^{k}(\g{so}\big(2n,\mathbb {R})^*\big)\tens \L\big(\g g_{1}^*\big)\Big)^{\g{sl}(2,\mathbb{R})}$ such that:
\begin{enumerate}
\item for $X\in \g{sl}(2,\mathbb{R})\setminus\NC$, $\phi(X)=\Spf(X)\in\sdir_{k=0}^n S^{k}(\g{so}\big(2n,\mathbb {R})^*\big)\tens \L\big(\g g_{1}^*\big)$;

\item $\W\phi=0$.
\end{enumerate}
Then we have $\phi=\Spf$.

\end{theo}

\end{document}